\documentclass[12pt]{article}
\usepackage{amsfonts,amsthm,amsmath,amssymb}
\usepackage{latexsym,graphicx,fullpage}
\usepackage{hypernat,titlesec,titletoc}

\title{{Minimum tree-stretch of Hamming \\
graphs and higher-dimensional grids}
\thanks{Supported by NSFC (61373106) and 973 Program of China (2010CB328101).}}
\author{{Lan Lin$^{1}$\thanks{Corresponding Author. E-mail address:
linlan@tongji.edu.cn} \ \ \ \ Yixun Lin$^{2}$}\\
{\small $^{1}$ School of Electronics and Information Engineering,
Tongji University,}\\ {\small Shanghai 200092, China}\\
{\small $^{2}$ School of Mathematics and Statistics, Zhengzhou University,}\\
{\small Zhengzhou 450001, China}}

\date{}
\begin{document}
\maketitle
\renewcommand\baselinestretch{1.0}

{\bf Abstract:} \ The minimum stretch spanning tree problem for a graph $G$
is to find a spanning tree $T$ of $G$ such that the maximum distance in $T$
between two adjacent vertices is minimized. The minimum value of this
optimization problem gives rise to a graph invariant $\sigma_T(G)$, called
the {\it tree-stretch} of $G$. The problem has been studied in the algorithmic
aspects, such as NP-hardness and fixed-parameter solvability. This paper
presents the exact values $\sigma_T(G)$ of the Hamming graphs $K_{n_1}\times
K_{n_2}\times \cdots \times K_{n_d}$ and the higher-dimensional grids
$P_{n_1}\times P_{n_2}\times \cdots \times P_{n_d}$.

{\bf Keywords:} \  spanning tree optimization, tree-stretch, tree-congestion,
Hamming graphs, higher-dimensional grids.

\section{Introduction }
\hspace*{0.3cm} Let $G$ be a simple connected graph with vertex set
$V(G)$ and edge set $E(G)$. Then $G$ contains a spanning tree $T$.
For $uv\in E(G)$, let $d_T(u,v)$ denote the distance between $u$ and
$v$ in $T$. The {\it max-stretch} of a spanning tree $T$ is defined by
\begin{equation}
\sigma_T(G,T):=\max_{uv\in E(G)}\,d_T(u,v).
\end{equation}
{\it The minimum stretch spanning tree problem} (the MSST problem
for short) is to find a spanning tree $T$ such that $\sigma_T(G,T)$
is minimized, where the minimum value
\begin{equation}
\sigma_T(G):=\min \{\sigma_T(G,T): T \mbox{ is a spanning tree of }
G\}
\end{equation}
is called the {\it tree-stretch} of $G$ (following the terminology of
\cite{Lieb08} and the notation $\sigma_T(G)$ of \cite{Fekete01}). A spanning
tree $T$ attaining this minimum value is called an {\it optimal spanning tree}.

With applications in distribution systems and communication networks, a series
of tree spanner problems were intensively studied in the literature (see
\cite{Peleg89,Cai95,Lieb08}, etc). A basic decision version of these problems
can be stated as follows: For a given integer $k$, is there a spanning tree $T$
of $G$ (called a tree $k$-spanner) such that the distance in $T$ between every
pair of vertices is at most $k$ times their distance in $G$? The MSST problem
mentioned above is the optimization version of this decision problem. Moreover,
it is worth pointing out that this graph embedding problem can be regarded as
a variant of the bandwidth (dilation) problem when graph $G$ is embedded into
its spanning tree $T$ (see survey \cite{Diaz02}).

In the algorithmic aspect, the MSST problem has been proved NP-hard (see
\cite{Brand04,Brand07,Cai95,Fekete01,Galb03}), and fixed-parameter
polynomial algorithms were discussed in details (see \cite{Brand04,Brand07,
Fekete01,Fomin11}). Several exact results of $\sigma_T(G)$ for special
graphs were also investigated. For example, the characterization of
$\sigma_T(G)=2$ was given in \cite{Bondy89,Cai95}. Besides, $\sigma_T(G)\leq 3$
for interval, split, and permutation graphs were showed in \cite{Brand04,
Le99,Madan96,Venka97}. Some formulas for basic families of special graphs,
such as complete $k$-partite graphs $K_{n_1,n_2,\ldots,n_k}$, rectangular grids
$P_m\times P_n$, torus grids $C_m\times C_n$, triangular grids $T_n$ and
hypercubes $Q_n$, can be seen in \cite{Lin17,Lin18}.

The MSST problem has close relations to {\it the minimum congestion spanning
tree problem} \cite{Ostrov04}, which is to find a spanning tree $T$ of $G$ such
that the size of the maximum fundamental edge-cut is minimized. The problem has
been proved NP-hard in \cite{Okamoto11}, and fixed-parameter polynomial algorithms
were presented in \cite{Bodlae12}. Much interest was paid to the exact
results for special graphs (see, e.g., \cite{Bodlae11,Caste09,Hruska08,Ostrov10}).
These results motivate our study on the MSST problem for typical graphs.

The embedding problems of Hamming graphs $K_{n_1}\times K_{n_2}\times\cdots
\times K_{n_d}$ have significant applications in error-correcting code and
multichannel communication. Even so, the bandwidth problem of Hamming graphs
was a long-standing open problem in this field (see \cite{Diaz02}). It is also
unsolved for the minimum congestion spanning tree problem. On the other hand,
various grid graphs are also appealing in graph embedding. For example, the minimum
tree-congestion and the minimum tree-stretch of two-dimensional grids $P_m\times P_n$
have been determined in \cite{Hruska08,Caste09,Lin17}, but the results for
higher-dimensional grids $P_{n_1}\times P_{n_2}\times\cdots \times P_{n_d}$ are
unknown yet. The goal of this paper is to determine the exact value $\sigma_T(G)$
for the Hamming graphs and the higher-dimensional grids.

The paper is organized as follows. In Section 2, some definitions and elementary
properties are introduced. In Section 3, we are concerned with the Hamming graphs.
Section 4 is devoted to the higher-dimensional grids. We give a short summary in
Section 5.

\section{Preliminaries}
\hspace*{0.3cm} We shall follow the graph-theoretic terminology and
notation of \cite{Bondy08}. Let $G$ be a simple connected graph on
$n$ vertices with vertex set $V(G)$ and edge set $E(G)$. For
$S\subseteq V(G)$, we denote by $G[S]$ the subgraph of $G$ induced
by $S$. For an edge $e\in E(G)$, denote by $G-e$ the graph obtained
from $G$ by deletion of $e$. For an edge $e$ not in $E(G)$, denote
by $G+e$ the graph obtained from $G$ by addition of $e$.

Let $P_n,C_n,K_n$ denote the path, the cycle, the complete graph,
respectively, on $n$ vertices. The {\it cartesian product} of two
graphs $G$ and $H$, denoted $G\times H$, is the graph with vertex set
$V(G)\times V(H)$ and two vertices $(u,v)$ and $(u',v')$ are adjacent
if and only if either $[u=u'\, \mbox{and}\, vv'\in E(H)]$ or $[v=v'\,
\mbox{and}\, uu'\in E(G)]$.

Let $T$ be a spanning tree of $G$. Usually, the spanning tree $T$ is
regarded as a set of edges. The {\it cotree} $\overline T$ of $T$ is the
complement of $T$ in $E(G)$, namely $\overline T=E(G)\setminus T$. For
an edge $e\in \overline T$, the unique cycle in $T+e$ is called the
{\it fundamental cycle} with respect to $e$. Moreover, the {\it detour}
for an edge $uv\in E(G)$ is the unique $u$-$v$ path in $T$, denoted by
$P_T(u,v)$. For $uv\in \overline T$, the fundamental cycle with respect
to $uv$ is indeed the detour $P_T(u,v)$ plus edge $uv$. So, the MSST
problem is equivalent to a problem of finding a spanning tree such that
the length of maximum fundamental cycle is minimized, where the tree-stretch
$\sigma_T(G)$ is one less than the length of this fundamental cycle.
As stated in \cite{Galb03}, since all fundamental cycles with respect
to a spanning tree $T$ constitute a basis of the {\it cycle space} of $G$,
this problem is indeed an optimal basis problem (which minimizes the
length of the longest cycle) in the cycle space.

For an edge $e\in T$, $T-e$ has two components; let $X_e$ be the vertex set
of one of these components. Then $\partial (X_e):=\{uv\in E(G):u\in X_e,
v\notin X_e\}$ is called the {\it fundamental edge-cut} (or {\it bond})
with respect to the tree-edge $e$. Here, $|\partial (X_e)|$ is called
the {\it congestion} of edge $e$. The minimum congestion spanning tree
problem is to find a spanning tree $T$ of $G$ such that the maximum
congestion in $T$ is minimized. This is an optimal basis problem in
the cocycle space.

The duality relation of the above problems lies on the following fact:
For $e'\in \overline T$, $e$ is contained in the fundamental cycle with
respect to $e'$ if and only if $e'$ is contained in the fundamental
edge-cut $\partial (X_e)$ with respect to $e$ (see \cite{Bondy08}).

The following observation is immediate.

{\bf Proposition 2.1} \ For a spanning tree $T$ of $G$, let $D(T)$ be the
diameter of $T$ (i.e., the maximum distance between any two vertices of $T$).
Then
$$ \sigma_T(G,T)\leq D(T).$$

{\bf Proof:} \ This is because for any $uv\in E(G)$, $d_T(u,v)\leq D(T)$. $\Box$

In a Hamming graph $G=K_{n_1}\times K_{n_2}\times \cdots \times K_{n_d}$,
each vertex can be represented by a $d$-dimensional vector $v=(x_1,x_2,\ldots,
x_d)$ with $x_i\in \{0,1,\ldots,n_i-1\}, n_i\geq 2$ for $1\leq i\leq d$. Two
vertices $v=(x_1,x_2,\ldots,x_d)$ and $u=(y_1,y_2,\ldots,y_d)$ are adjacent if
they differ in exactly one coordinate. An illustration is shown in Figure 1.

\begin{center}
\setlength{\unitlength}{0.32cm}
\begin{picture}(31,15)

\multiput(1,4)(3,0){5}{\circle*{0.3}}
\multiput(1,7)(3,0){5}{\circle*{0.3}}
\multiput(1,10)(3,0){5}{\circle*{0.3}}
\multiput(1,13)(3,0){5}{\circle*{0.3}}

\put(1,4){\line(1,0){12}} \put(1,7){\line(1,0){12}}
\put(1,10){\line(1,0){12}}\put(1,13){\line(1,0){12}}
\put(1,4){\line(0,1){9}} \put(4,4){\line(0,1){9}}
\put(7,4){\line(0,1){9}} \put(10,4){\line(0,1){9}}
\put(13,4){\line(0,1){9}}
\bezier{40}(1,13)(7,15.5)(13,13)
\bezier{30}(1,13)(6.5,14.4)(10,13)
\bezier{30}(4,13)(8.5,14.5)(13,13)
\bezier{30}(1,13)(4,14)(7,13)
\bezier{30}(4,13)(7,14)(10,13)
\bezier{30}(7,13)(10,14)(13,13)
\bezier{40}(1,10)(7,12.5)(13,10)
\bezier{30}(1,10)(6.5,11.4)(10,10)
\bezier{30}(4,10)(8.5,11.5)(13,10)
\bezier{30}(1,10)(4,11)(7,10)
\bezier{30}(4,10)(7,11)(10,10)
\bezier{30}(7,10)(10,11)(13,10)
\bezier{40}(1,7)(7,9.5)(13,7)
\bezier{30}(1,7)(6.5,8.4)(10,7)
\bezier{30}(4,7)(8.5,8.5)(13,7)
\bezier{30}(1,7)(4,8)(7,7)
\bezier{30}(4,7)(7,8)(10,7)
\bezier{30}(7,7)(10,8)(13,7)
\bezier{40}(1,4)(7,6.5)(13,4)
\bezier{30}(1,4)(6.5,5.4)(10,4)
\bezier{30}(4,4)(8.5,5.5)(13,4)
\bezier{30}(1,4)(4,5)(7,4)
\bezier{30}(4,4)(7,5)(10,4)
\bezier{30}(7,4)(10,5)(13,4)

\bezier{30}(1,4)(-1,8.5)(1,13)
\bezier{30}(1,4)(-0.2,7)(1,10)
\bezier{30}(1,7)(-0.2,10)(1,13)
\bezier{30}(4,4)(2,8.5)(4,13)
\bezier{30}(4,4)(2.8,7)(4,10)
\bezier{30}(4,7)(2.8,10)(4,13)
\bezier{30}(7,4)(5,8.5)(7,13)
\bezier{30}(7,4)(5.8,7)(7,10)
\bezier{30}(7,7)(5.8,10)(7,13)
\bezier{30}(10,4)(8,8.5)(10,13)
\bezier{30}(10,4)(8.8,7)(10,10)
\bezier{30}(10,7)(8.8,10)(10,13)
\bezier{30}(13,4)(11,8.5)(13,13)
\bezier{30}(13,4)(11.8,7)(13,10)
\bezier{30}(13,7)(11.8,10)(13,13)

\put(19,4){\line(1,0){8}} \put(19,7){\line(1,0){8}}
\put(19,10){\line(1,0){8}}\put(19,4){\line(1,1){4}}
\put(19,7){\line(1,1){4}} \put(19,10){\line(1,1){4}}
\put(23,8){\line(1,0){8}} \put(23,11){\line(1,0){8}}
\put(23,14){\line(1,0){8}}\put(27,4){\line(1,1){4}}
\put(27,7){\line(1,1){4}} \put(27,10){\line(1,1){4}}
\put(19,4){\line(0,1){6}} \put(23,8){\line(0,1){6}}
\put(27,4){\line(0,1){6}}\put(31,8){\line(0,1){6}}

\put(-1,8){\makebox(1,0.5)[l]{\small $K_4$}}
\put(6,14){\makebox(1,0.5)[l]{\small $K_5$}}
\put(17.3,7){\makebox(1,0.5)[l]{\small $K_l$}}
\put(19.2,11.8){\makebox(1,0.5)[l]{\small $K_m$}}
\put(26.5,14.5){\makebox(1,0.5)[l]{\small $K_n$}}
\put(22,5.5){\makebox(1,0.5)[l]{\small $K_m\times K_n$}}

\put(3,1.8){\makebox(1,0.5)[l]{\small (a) $K_4\times K_5$}}
\put(20,1.8){\makebox(1,0.5)[l]{\small (b) $K_l\times K_m\times K_n$}}

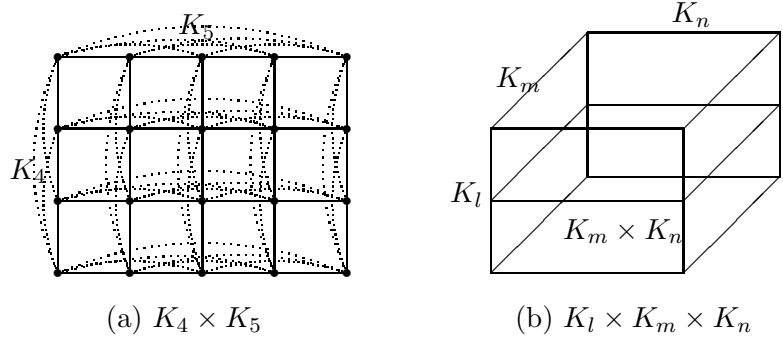
\put(4,0){\makebox(1,0.5)[l]{\small Figure 1. Illustration of Hamming graphs }}
\end{picture}
\end{center}

When $n_1=n_2=\cdots=n_d=2$, $Q_d=K_2\times K_2\times \cdots \times K_2$ is called
a $d$-dimensional {\it hypercube} or $d$-{\it cube}, whose vertex set is the set
of all $d$-dimensional 0-1 vectors $(x_1,x_2,\ldots,x_d)$ and two vertices are
adjacent if they differ in exactly one coordinate.

In a $d$-dimensional grid $P_{n_1}\times P_{n_2}\times\cdots \times P_{n_d}$,
each vertex can be represented by a $d$-dimensional vector $v=(x_1,x_2,\ldots,
x_d)$ with $x_i\in \{1,2,\ldots,n_i\}, n_i\geq 2$ for $1\leq i\leq d$. Two
vertices $v=(x_1,x_2,\ldots,x_d)$ and $u=(y_1,y_2,\ldots,y_d)$ are adjacent if
they differ by 1 in exactly one coordinate. The illustration of higher-dimensional
grids is similar to Figure 1 with $K_{n_i}$ replaced by $P_{n_i}$ for $1\leq i\leq d$.

\section{Hamming graphs}
\hspace*{0.3cm} On the one hand, we shall show the upper bound of
$\sigma_T(K_{n_1}\times K_{n_2}\times \cdots \times K_{n_d})$.

{\bf Lemma 3.1} \ For the Hamming graphs $G=K_{n_1}\times K_{n_2}
\times\cdots \times K_{n_d}$ ($2\leq n_1\leq n_2\leq \cdots \leq
n_d$), it holds that
$$\sigma_T(K_{n_1}\times K_{n_2}\times\cdots \times K_{n_d})
\leq \begin{cases}
2d-1,& \mbox{if}\,\, n_1=2\\
2d,& \mbox{if}\,\, n_1\geq 3.
\end{cases}$$

{\bf Proof:} \ We first consider $d=2$ and $G=K_{n_1}\times K_{n_2}$. Suppose
that $V(G):=\{(i,j):0\leq i\leq n_1-1,0\leq j\leq n_2-1\}$ and $(i,j)$ is
adjacent to $(i',j')$ if $i=i'$ or $j=j'$. We may call $R_i:=\{(i,j):
0\leq j\leq n_2-1\}$ the $i$-th row for $0\leq i\leq n_1-1$, and $C_j:=\{(i,j):
0\leq i\leq n_1-1\}$ the $j$-th column for $0\leq j\leq n_2-1$, each of which
is a clique. We construct a spanning tree $T_2$ as follows. First, take a star
(regarded as $T_1$) in each row $R_i$ with the center at column $C_1$ ($0\leq
i\leq n_1-1$). Then, take a star in column $C_1$ with center $x_0=(0,0)$ to
join the centers of stars in rows. An example is shown in Figure 2. For $n_1=2$,
$T_2$ is a tree of diameter three (double star), and so by Proposition 2.1, we
have $\sigma_T(G)\leq\sigma_T(G,T_2)\leq 3$. For $n_1\geq 3$, $T_2$ is a tree
of diameter four (each leaf has distance at most two from $x_0$). By Proposition 2.1,
it follows that $\sigma_T(G)\leq \sigma_T(G,T_2)\leq 4$. Therefore, by means
of this spanning tree $T_2$, we have
$$\sigma_T(G)\leq \begin{cases}
3,& \mbox{if}\,\, n_1=2\\
4,& \mbox{if}\,\, n_1\geq 3,
\end{cases}$$
and so the assertion is true for $d=2$.

\begin{center}
\setlength{\unitlength}{0.32cm}
\begin{picture}(14,15)

\multiput(1,4)(3,0){5}{\circle*{0.3}}
\multiput(1,7)(3,0){5}{\circle*{0.3}}
\multiput(1,10)(3,0){5}{\circle*{0.3}}
\multiput(1,13)(3,0){5}{\circle*{0.3}}

\put(1,4){\line(1,0){3}} \put(1,7){\line(1,0){3}}
\put(1,10){\line(1,0){3}}\put(1,13){\line(1,0){3}}
\put(1,10){\line(0,1){3}}

\qbezier(1,4)(4,5)(7,4) \qbezier(1,4)(5.5,6)(10,4)
\qbezier(1,4)(7,7)(13,4) \qbezier(1,7)(4,8)(7,7)
\qbezier(1,7)(5.5,9)(10,7) \qbezier(1,7)(7,10)(13,7)
\qbezier(1,10)(4,11)(7,10) \qbezier(1,10)(5.5,12)(10,10)
\qbezier(1,10)(7,13)(13,10) \qbezier(1,13)(4,14)(7,13)
\qbezier(1,13)(5.5,15)(10,13) \qbezier(1,13)(7,16)(13,13)
\qbezier(1,4)(-1,8.5)(1,13) \qbezier(1,7)(0,10)(1,13)

\put(-0.5,13){\makebox(1,0.5)[l]{\small $x_0$}}
\put(-2,9){\makebox(1,0.5)[l]{\small $K_4$}}
\put(7,15){\makebox(1,0.5)[l]{\small $K_5$}}

\put(-3,0){\makebox(1,0.5)[l]{\small Figure 2. Spanning tree for
$K_4\times K_5$}}
\end{picture}
\end{center}

We proceed to construct a spanning tree $T_d$ by induction on $d$.
For the spanning tree $T_2$ before, we call the vertex $x_0$ the
{\it center} of $T_2$ and denoted $x^2_0$ henceforth. Assume that
$d\geq 3$ and $T_{d-1}$ has been constructed. The graph $G$ now
consists of $n_1$ copies of $K_{n_2}\times\cdots \times K_{n_d}$,
each of which has a required spanning tree $T_{d-1}$. We construct
a spanning tree $T_d$ of $G$ by joining a star between $n_1$ centers
of the copies $T_{d-1}$; and let $x^d_0$ be the center of this star.
Then $x^d_0$ is called the center of $T_d$.

This spanning tree $T_d$ has the property that each leaf of
$T_d$ has distance at most $d$ from the center $x^d_0$. In fact, this
is trivially true for $d=1,2$ (see Figure 2). If it is true for
$T_{d-1}$, that is, each leaf has distance at most $d-1$ from the
center $x^{d-1}_0$, then, since $x^{d-1}_0$ and $x^d_0$ have distance
one, the property follows for $T_d$.

To show the upper bound, we first consider the case $n_1=2$. There are
only two copies of $T_{d-1}$ in $G$. Then $T_d$ is obtained by joining
an edge between the centers of these two copies of $T_{d-1}$. Since each
leaf of $T_{d-1}$ has distance at most $d-1$ from the center, it follows
that the diameter of $T_d$ is at most $2(d-1)+1=2d-1$. By Proposition 2.1,
we see that $\sigma_T(G)\leq \sigma_T(G,T_d)\leq D(T_d)\leq 2d-1$.

We next consider the case $n_1\geq 3$. Now, $T_d$ is constructed as follows:
Among the $n_1$ centers of copies $T_{d-1}$, we choose one as the center
$x^d_0$ of $T_d$, and join a star connecting to the other $n_1-1$ centers
of copies $T_{d-1}$ (as the star in column $C_1$ of Figure 2).  Since each
leaf of $T_{d-1}$ has distance at most $d-1$ from the center, it follows
that the diameter of $T_d$ is at most $2(d-1)+2=2d$. By Proposition 2.1,
we have $\sigma_T(G)\leq \sigma_T(G,T_d)\leq D(T_d)\leq 2d$. Thus the
assertion is proved.  $\Box$

On the other hand, we shall show the lower bound. In the Hamming
graph $G=K_{n_1}\times K_{n_2}\times\cdots \times K_{n_d}$ ($n_i\geq 2$),
for each vertex $v=(x_1,x_2,\ldots,x_d)$, the vertex $f(v)=(x'_1,x'_2,
\ldots,x'_d)$ is called the {\it antipodal vertex} of $v$ if
\begin{equation}
x'_i=x_i+1\,(\mbox{ mod}\,n_i), \quad i=1,2,\ldots,d.
\end{equation}

Obviously, $v$ and $f(v)$ have distance $d$ in $G$. Note that
this definition is not symmetrical, as the antipodal vertex of $f(v)$
is not necessarily $v$.

{\bf Theorem 3.2} \ For the Hamming graphs $G=K_{n_1}\times K_{n_2}
\times \cdots \times K_{n_d}$ ($2\leq n_1\leq n_2\leq \cdots \leq
n_k$), it holds that
$$\sigma_T(K_{n_1}\times K_{n_2}\times\cdots \times K_{n_d})
= \begin{cases}
2d-1,& \mbox{if}\,\, n_1=2\\
2d,& \mbox{if}\,\, n_1\geq 3.
\end{cases}$$

{\bf Proof:} \ By Lemma 3.1, it suffices to show the lower bound. Suppose that
$T$ is an arbitrary spanning tree of $G$. We distinguish two cases as follows.

{\bf Case 1:} $n_1=2$. We shall show that $\sigma_T(G,T)\geq 2d-1$.

For each vertex $v\in V(G)$, let $f(v)$ be the antipodal vertex of $v$.
We denote by $P_T(v,f(v))$ the path in $T$ from $v$ to $f(v)$.
Furthermore, we define the {\it successor} of $v$, denoted $s(v)$, by the
next vertex of $v$ on the path $P_T(v,f(v))$. We claim that there exists an
edge $uv\in T$ such that $s(u)=v$ and $s(v)=u$. Assume, to the contrary,
that there is no such edge. Then we can start at a vertex $v_1$ and let
$v_2=s(v_1)$. Since $s(v_2)\neq v_1$, $v_3=s(v_2)$ is a new vertex. In
this way, we can define a sequence $(v_1,v_2,v_3,\ldots)$ by setting
$v_{i+1}=s(v_i)$ in the spanning tree $T$. Since $T$ contains no cycles,
each vertex in this sequence cannot repeat the ones previously visited.
Hence this is indeed an infinite sequence, contradicting that $T$ is a
finite tree.

Now we take an edge $uv\in T$ that $s(u)=v$ and $s(v)=u$. Then $f(u)$ and
$f(v)$ belong to different components of $T-uv$. Thus the paths $P_T(u,f(u))$
and $P_T(v,f(v))$ have only the edge $uv$ in common. We take $P_T(f(u),f(v))
=P_T(f(u),u)\cup P_T(v,f(v))$ by joining these two paths. Suppose that
$$u=(x_1,\ldots,x_{i-1},a,x_{i+1},\ldots,x_d),$$
$$v=(x_1,\ldots,x_{i-1},b,x_{i+1},\ldots,x_d),$$
where $a,b\in \{0,1,\ldots,n_i-1\}$, and $a\neq b$. Then
$$f(u)=(x_1+1,\ldots,x_{i-1}+1,a+1,x_{i+1}+1,\ldots,x_d+1),$$
$$f(v)=(x_1+1,\ldots,x_{i-1}+1,b+1,x_{i+1}+1,\ldots,x_d+1),$$
where the addition is modular as in (3). So $f(u)$ and $f(v)$ are adjacent
in $G$. Furthermore, $uv\in T$ implies $f(u)f(v)\in \overline T$ (for otherwise
there would be a cycle in $T$). Note that the lengths of $P_T(u,f(u))$ and
$P_T(v,f(v))$ are at least $d$, and they have $uv$ in common. Hence the
length of $P_T(f(u),f(v))$ is at least $2d-1$, and so $d_T(f(u),f(v))\geq 2d-1$.
Therefore, we deduce the lower bound $\sigma_T(G,T)\geq 2d-1$.

{\bf Case 2:} $n_1\geq 3$. We shall show that $\sigma_T(G,T)\geq 2d$.

By the proof of Case 1, we found a path $P_T(f(u),f(v))$ with $d_T(f(u),f(v))
\geq 2d-1$. If $d_T(f(u),f(v))>2d-1$, then we are done. So we may assume that
$d_T(f(u),f(v))=2d-1$. Thus $d_T(u,f(u))=d_T(v,f(v))=d$. We observe the path
$P_T(u,f(u))$ from $u=(x_1,\ldots,x_{i-1},a,x_{i+1},\ldots,x_d)$ to $f(u)=
(x_1+1,\ldots,x_{i-1}+1,a+1,x_{i+1}+1,\ldots,x_d+1)$. Since its length is exactly
$d$, it follows that each coordinate of the $d$-dimensional vector increases
exactly by 1 (in modular sense) successively along this path. Noting that $v=s(u)$
is the next vertex of $u$ on this path, we see that $b=a+1\,(\mbox{ mod}\,n_i)$.
However, as $n_i\geq n_1\geq 3$, we know that $b+1\neq a\,(\mbox{ mod}\,n_i)$.

Let us now observe another path $P_T(v,f(v))$ from $v=(x_1,\ldots,x_{i-1},b,x_{i+1},
\ldots,x_d)$ to $f(v)=(x_1+1,\ldots,x_{i-1}+1,b+1,x_{i+1}+1,\ldots,x_d+1)$. Here, the
next vertex of $v$ on this path is $u$, instead of $(x_1,\ldots,x_{i-1},b+1,x_{i+1},
\ldots,x_d)$ (note that $a\neq b+1$). Along this path after vertex $u$, the coordinate
$x_i=a$ of $u$ must have one more changing to $x_i=b+1$ somewhere. Consequently,
the length of $P_T(v,f(v))$ is greater than $d$. That is to say, the assumption of
$d_T(v,f(v))=d$ is impossible. Therefore, $d_T(f(u),f(v))>2d-1$, which gives the
lower bound $\sigma_T(G,T)\geq 2d$.

Combining the lower bound here and the upper bound in Lemma 3.1 completes the proof
of the theorem. $\Box$

As a special case, we derive the following result for hypercubes in \cite{Lin18}. In the
above proof, we also generalize a property of R.L. Graham on hypercubes (every spanning tree
of $Q_d$ has a fundamental cycle of length at least $2d$, see Exercise 4.2.15(d) of
\cite{Bondy08}) to the Hamming graphs.

{\bf Corollary 3.3} \ For the hypercubes $Q_d$, it holds that $\sigma_T(Q_d)=2d-1$.

\section{Higher-dimensional grids}
\hspace*{0.3cm} We consider the $d$-dimensional grids $P_{n_1}\times P_{n_2}\times\cdots \times
P_{n_d}$, in which each vertex is represented by $v=(x_1,x_2,\ldots,x_d)$ with
$x_i\in \{1,2,\ldots,n_i\}, n_i\geq 2$ for $1\leq i\leq d$. Two vertices
$v=(x_1,x_2,\ldots,x_d)$ and $u=(y_1,y_2,\ldots,y_d)$ are adjacent if they differ
by 1 in exactly one coordinate.

Suppose that $G_i=P_{n_1}\times P_{n_2}\times\cdots \times P_{n_i}$ ($1\leq i\leq d$).
Then $G_1=P_{n_1}$, $G_i=G_{i-1}\times P_{n_i}$ ($i\geq 2$), and $G_d=G$. When $d=1$,
it is trivial that $G_1=P_{n_1}$ is a path and $T_1=G_1$ is an optimal tree. In the
sequel, for a path $P_n=(v_1,v_2,\ldots,v_n)$, the vertex $v_{\lfloor n/2\rfloor}$ is
called the {\it center} of this path.

We begin with the case $d=2$. For a $2$-dimensional grid $G=P_{n_1}\times P_{n_2}$
($2\leq n_1\leq n_2$), let $V(G):=\{(i,j):1\leq i\leq n_1,1\leq j\leq n_2\}$
denote the vertex set, and $(i,j)$ is adjacent to $(i',j')$ if $|i-i'|+|j-j|=1$.
We call $R_i:=\{(i,j):1\leq j\leq n_2\}$ the $i$-th row, and $C_j:=\{(i,j):1\leq
i\leq n_1\}$ the $j$-th column.

We construct a spanning tree $T_2$ as follows. First, we take $n_2$ copies of $P_{n_1}$,
namely, the $n_2$ columns. Let $v_1,v_2,\ldots,v_{n_2}$ be the centers of these columns
in turn. Then we join a path $P_{n_2}$ to pass through these $n_2$ centers, namely,
the row $R_{\lfloor n_1/2\rfloor}$. An example is shown in Figure 3 (where the centers
are marked by heavy dots). In this tree $T_2$, the path $P_{n_2}$ passing through $n_2$
centers of $P_{n_1}$'s is called the {\it central path} (see $R_2$ in Figure 3). The
center of this central path is called the {\it center} of $T_2$.

\begin{center}
\setlength{\unitlength}{0.32cm}
\begin{picture}(15,14)

\multiput(1,4)(3,0){5}{\circle*{0.3}}
\multiput(1,7)(3,0){5}{\circle*{0.3}}
\multiput(1,10)(3,0){5}{\circle*{0.3}}
\multiput(1,10)(3,0){5}{\circle{0.4}}
\multiput(7,10)(0,0){1}{\circle{0.5}}
\multiput(1,13)(3,0){5}{\circle*{0.3}}

\put(1,10){\line(1,0){12}}
\put(1,4){\line(0,1){9}} \put(4,4){\line(0,1){9}}
\put(7,4){\line(0,1){9}} \put(10,4){\line(0,1){9}}
\put(13,4){\line(0,1){9}}
\bezier{30}(1,4)(7,4)(13,4)
\bezier{30}(1,7)(7,7)(13,7)
\bezier{30}(1,13)(7,13)(13,13)

\put(-1,12.8){\makebox(1,0.5)[l]{\small $R_1$}}
\put(-1,9.8){\makebox(1,0.5)[l]{\small $R_2$}}
\put(-1,6.8){\makebox(1,0.5)[l]{\small $R_3$}}
\put(-1,3.8){\makebox(1,0.5)[l]{\small $R_4$}}
\put(0.6,2.5){\makebox(1,0.5)[l]{\small $C_1$}}
\put(3.6,2.5){\makebox(1,0.5)[l]{\small $C_2$}}
\put(6.6,2.5){\makebox(1,0.5)[l]{\small $C_3$}}
\put(9.6,2.5){\makebox(1,0.5)[l]{\small $C_4$}}
\put(12.6,2.5){\makebox(1,0.5)[l]{\small $C_5$}}
\put(-2,0){\makebox(1,0.5)[l]{\small
Figure 3. Spanning tree for $P_4\times P_5$.}}
\end{picture}
\end{center}

As we shall see later in the general case, this spanning tree $T_2$ is optimal, and a detour
in $T_2$ with the maximum stretch is as follows: Start at a leaf of column $C_i$ with
distance $\left\lfloor\frac{n_1}{2}\right\rfloor$ to the central path, go to the central
path, pass through an edge of it, and then go back to the leaf of column $C_{i+1}$ at the
same side. So we have the following (see \cite{Lin17}):
$$\sigma_T(P_{n_1}\times P_{n_2})=2\left\lfloor\frac{n_1}{2}\right\rfloor+1.$$

\begin{center}
\setlength{\unitlength}{0.32cm}
\begin{picture}(22,20)

\multiput(0,4)(6,0){4}{\circle*{0.3}}
\multiput(2,6)(6,0){4}{\circle*{0.3}}
\multiput(4,8)(6,0){4}{\circle*{0.3}}
\multiput(0,8)(6,0){4}{\circle*{0.3}}
\multiput(2,10)(6,0){4}{\circle*{0.3}}
\multiput(4,12)(6,0){4}{\circle*{0.3}}
\multiput(0,12)(6,0){4}{\circle*{0.3}}
\multiput(2,14)(6,0){4}{\circle*{0.3}}
\multiput(2,14)(6,0){4}{\circle{0.4}}
\multiput(8,14)(0,0){1}{\circle{0.5}}
\multiput(4,16)(6,0){4}{\circle*{0.3}}
\multiput(0,16)(6,0){4}{\circle*{0.3}}
\multiput(2,18)(6,0){4}{\circle*{0.3}}
\multiput(4,20)(6,0){4}{\circle*{0.3}}

\put(0,4){\line(1,1){4}} \put(0,8){\line(1,1){4}}
\put(0,12){\line(1,1){4}} \put(0,16){\line(1,1){4}}
\put(2,6){\line(0,1){12}}

\put(6,4){\line(1,1){4}} \put(6,8){\line(1,1){4}}
\put(6,12){\line(1,1){4}} \put(6,16){\line(1,1){4}}
\put(8,6){\line(0,1){12}}

\put(12,4){\line(1,1){4}} \put(12,8){\line(1,1){4}}
\put(12,12){\line(1,1){4}} \put(12,16){\line(1,1){4}}
\put(14,6){\line(0,1){12}}

\put(18,4){\line(1,1){4}} \put(18,8){\line(1,1){4}}
\put(18,12){\line(1,1){4}} \put(18,16){\line(1,1){4}}
\put(20,6){\line(0,1){12}}

\put(2,14){\line(1,0){18}}
\bezier{50}(0,4)(9,4)(18,4)
\bezier{50}(4,8)(13,8)(22,8)
\bezier{50}(0,16)(9,16)(18,16)
\bezier{50}(4,20)(13,20)(22,20)

\bezier{50}(0,4)(0,10)(0,16)
\bezier{50}(6,4)(6,10)(6,16)
\bezier{50}(12,4)(12,10)(12,16)
\bezier{50}(18,4)(18,10)(18,16)

\bezier{50}(4,8)(4,14)(4,20)
\bezier{50}(10,8)(10,14)(10,20)
\bezier{50}(16,8)(16,14)(16,20)
\bezier{50}(22,8)(22,14)(22,20)

\put(0,0){\makebox(1,0.5)[l]{\small
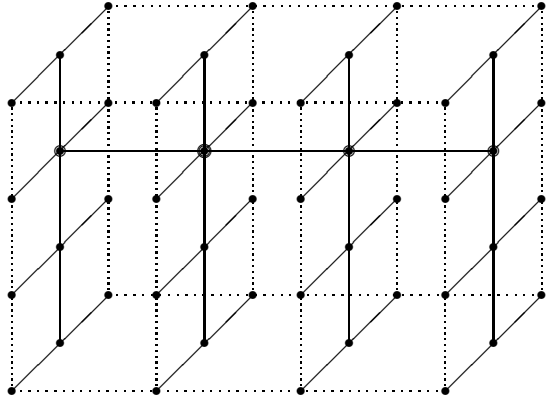
Figure 4. Spanning tree for $P_3\times P_4\times P_4$.}}
\end{picture}
\end{center}

We further consider the case $d\geq 3$. For this, a spanning tree $T_i$ of $G_i$ in constructed
by induction on $i$. Assume that $i\geq 3$ and $T_{i-1}$ has been constructed. The graph $G_i$
now consists of $n_i$ copies of $G_{i-1}$, denoted $G^1_{i-1},G^2_{i-1},\ldots,G^{n_i}_{i-1}$,
where $G^j_{i-1}$ has a spanning tree $T^j_{i-1}$ ($1\leq j\leq n_i$). We construct a spanning
tree $T_i$ of $G_i$ by joining a central path $P_{n_i}$ through $n_i$ centers of $T^1_{i-1},
T^2_{i-1},\ldots,T^{n_i}_{i-1}$. The center of the central path $P_{n_i}$ is the center of
$T_i$, denoted $v^0_i$.  An example of $d=3$ is shown in Figure 4 (where the centers are marked
by heavy dots).

We proceed to show that $T_d$ is an optimal spanning tree of $G$. First, the
following lemma gives the upper bound.

{\bf Lemma 4.1} \ For the $d$-dimensional grids $G=P_{n_1}\times P_{n_2}
\times\cdots \times P_{n_d}$ ($2\leq n_1\leq n_2\leq \cdots \leq
n_d$), it holds that
$$\sigma_T(K_{n_1}\times K_{n_2}\times\cdots \times K_{n_d})
\leq 2\left(\left\lfloor\frac{n_1}{2}\right\rfloor+\left\lfloor\frac{n_2}{2}\right\rfloor
+\cdots+\left\lfloor\frac{n_{d-1}}{2}\right\rfloor \right)+1.$$

{\bf Proof:} \ We have defined the spanning tree $T_i$ of $G_i$ with center $v^0_i$
inductively. The following property plays an important role.

{\bf Claim} \ In the spanning tree $T_i$ of $G_i$, the distance between each
vertex $v$ and the center $v^0_i$ of $T_i$ satisfies
$$d_{T_i}(v,v^0_i)\leq \left\lfloor\frac{n_1}{2}\right\rfloor+\left\lfloor\frac{n_2}{2}\right\rfloor
+\cdots+\left\lfloor\frac{n_i}{2}\right\rfloor.$$

To see this, we use induction on $i$. When $i=1$, $T_1$ is a path and the distance
between each vertex and the center is at most $\lfloor n_1/2\rfloor$. Suppose
that $i\geq 2$ and the assertion is true for smaller $i$. Note that $G_i=
G_{i-1}\times P_{n_i}$ and $G_i$ consists of $n_i$ copies of $G_{i-1}$, that is,
$G^1_{i-1},G^2_{i-1},\ldots,G^{n_i}_{i-1}$, where $G^j_{i-1}$ has a spanning tree
$T^j_{i-1}$ ($1\leq j\leq n_i$). Let $v$ be a vertex of $G_i$. Then $v$ belongs to
some copy $G^j_{i-1}$. It follows that $P_{T_i}(v,v^0_i)=P_{T_i}(v,v^0_{i-1})\cup
P_{T_i}(v^0_{i-1},v^0_i)$, where $P_{T_i}(v,v^0_{i-1})=P_{T_{i-1}}(v,v^0_{i-1})$
and $P_{T_i}(v^0_{i-1},v^0_i)$ is contained in the central path $P_{n_i}$.  By
the inductive hypothesis, we have
$$d_{T_i}(v,v^0_{i-1})\leq \left\lfloor\frac{n_1}{2}\right\rfloor+\left\lfloor\frac{n_2}{2}\right\rfloor
+\cdots+\left\lfloor\frac{n_{i-1}}{2}\right\rfloor.$$
and
$$d_{T_i}(v^0_{i-1},v^0_i)\leq \left\lfloor\frac{n_i}{2}\right\rfloor.$$
Combining the above two inequalities results in the claim.

Now consider the spanning tree $T_d$ of $G$. For any cotree-edge $uv\in \overline T_d$,
we may assume that $uv$ is between two copies of $G_{d-1}$. For otherwise $uv$ is
contained in a copy of $G_{d-1}$, and we can get a smaller upper bound by the same
method. Suppose that $u^0_{d-1}$ and $v^0_{d-1}$ are the centers of these two copies
of $G_{d-1}$. Then $P_{T_d}(u,v)=P_{T_{d-1}}(u,u^0_{d-1})\cup \{u^0_{d-1}v^0_{d-1}\}
\cup P_{T_{d-1}}(v,v^0_{d-1})$, where the edge $u^0_{d-1}v^0_{d-1}$ is contained in
the central path. By the above Claim, we have
$$d_{T_d}(u,v)\leq 2\left(\left\lfloor\frac{n_1}{2}\right\rfloor+\left\lfloor\frac{n_2}{2}\right\rfloor
+\cdots+\left\lfloor\frac{n_{d-1}}{2}\right\rfloor\right)+1.$$
Since this inequality holds for every cotree-edge $uv\in \overline T_d$, we deduce
the upper bound in the theorem. $\Box$

It remains to show the lower bound. To this end, we pay attention to some special paths as
follows. For $a_j\in \{1,n_j\}, 1\leq j\leq d,j\neq i$, the path
$$P^A_{n_i}:=((a_1,\ldots,a_{i-1},1,a_{i+1},\ldots,a_d),(a_1,\ldots,a_{i-1},2,a_{i+1},\ldots,a_d),$$
$$\ldots,(a_1,\ldots,a_{i-1},n_i,a_{i+1},\ldots,a_d))$$
is called a {\it boundary path in $x_i$-coordinate}. For all possible combinations of
$a_j\in \{1,n_j\}, 1\leq j\leq d,j\neq i$, there are $2^{d-1}$ boundary paths in
$x_i$-coordinate. So there are totally $d2^{d-1}$ boundary paths in all coordinates.
For example, in the case $d=2$ of Figure 3, the 4 boundary paths are those on the 4
sides of the rectangular region. In the case $d=3$ of Figure 4, the 4 boundary paths
in $x_3$-coordinate are those on the horizontal dotted lines and all 12 boundary paths
are those on the 12 edges of the parallelepiped.

Moreover, the {\it antipodal boundary path} of $P^A_{n_i}$ is defined by
$$P^B_{n_i}:=((b_1,\ldots,b_{i-1},1,b_{i+1},\ldots,b_d),(b_1,\ldots,b_{i-1},2,b_{i+1},\ldots,b_d),$$
$$\ldots,(b_1,\ldots,b_{i-1},n_i,b_{i+1},\ldots,b_d))$$
where $\{a_j,b_j\}=\{1,n_j\}$ for $1\leq j\leq d,j\neq i$. For example, in the case
$d=2$ of Figure 3, two antipodal boundary paths are on the opposite sides of the
rectangular region. In the case $d=3$ of Figure 4, two antipodal boundary paths are
on the opposite edges of the parallelepiped (two opposite edges of a polyhedron are
two edges which are not contained in the same face).

By virtue of symmetry, the following antipodal boundary paths in $x_d$-coordinate 
are said to be {\it in standard form}:
\begin{eqnarray}
P^A_{n_d} &:=& ((1,1,\ldots,1,1),(1,1,\ldots,1,2),\ldots,(1,1,\ldots,1,n_d)),\\
P^B_{n_d} &:=& ((n_1,n_2,\ldots,n_{d-1},1),(n_1,n_2,\ldots,n_{d-1},2),\ldots,(n_1,n_2,\ldots,n_{d-1},n_d)).
\end{eqnarray}
In fact, any pair of antipodal boundary paths $P^A_{n_d}$ and $P^B_{n_d}$ can be
transformed into this form by reversing the order of some $x_i$-coordinates from
$(1,2,\ldots,n_{i-1},n_i)$ to $(n_i,n_{i-1},\ldots,2,1)$ if necessary.

In the proof below, we mainly analyze the relationship of two antipodal boundary paths.

{\bf Theorem 4.2} \ For the $d$-dimensional grids $G=P_{n_1}\times P_{n_2}
\times\cdots \times P_{n_d}$ ($2\leq n_1\leq n_2\leq \cdots \leq
n_d$), it holds that
$$\sigma_T(K_{n_1}\times K_{n_2}\times\cdots \times K_{n_d})
= 2\left(\left\lfloor\frac{n_1}{2}\right\rfloor+\left\lfloor\frac{n_2}{2}\right\rfloor
+\cdots+\left\lfloor\frac{n_{d-1}}{2}\right\rfloor\right)+1.$$

{\bf Proof:} \ By Lemma 4.1, it suffices to show the upper bound is also the lower bound.
Suppose that $T$ is an arbitrary spanning tree of $G$. Note that $2\leq n_1\leq n_2\leq \cdots \leq
n_d$. We consider $n_d$ copies $G^1_{d-1},G^2_{d-1},\ldots,G^{n_d}_{d-1}$ of $G=G_d$.
For clarity, we may denote these $(d-1)$-dimensional grids by $H_1=G^1_{d-1},H_2=G^2_{d-1},
\ldots,H_{n_d}=G^{n_d}_{d-1}$. Moreover, let $P_T(v_1,v_k)=(v_1,v_2,\ldots,v_k)$ be a shortest
path in $T$ from $H_1$ to $H_{n_d}$, where $v_1\in V(H_1)$, $v_k\in V(H_{n_d})$, and
$k\geq n_d$.

We take a pair of antipodal boundary paths $P^A_{n_d}$ and $P^B_{n_d}$ in $x_d$-coordinate as shown in
(4) and (5). They are also from $H_1$ to $H_{n_d}$, but not necessarily contained in $T$. We distinguish
three cases as follows.

{\bf Case 1:} \ Neither $P^A_{n_d}$ nor $P^B_{n_d}$ is contained in $T$.

For an edge $v_jv_{j+1}\in T$ in the path $P_T(v_1,v_k)$, let $S_j$ and $\overline S_j$
be the vertex sets of the two components of $T-v_jv_{j+1}$, where $v_j\in S_j$ and
$v_{j+1}\in \overline S_j$. Then $\partial (S_j):=\{uv\in E(G):u\in S_j,v\in \overline S_j\}$
is a fundamental edge-cut with respect to the tree-edge $v_jv_{j+1}\in T$. Let $u_1=(1,1,
\ldots,1)$ be the first vertex of $P^A_{n_d}$ and $w_1=(n_1,n_2,\ldots,n_{d-1},1)$ the
first vertex of $P^B_{n_d}$. Suppose that $P_T(u_1,v_{i_1})$ is the path in $T$ from
$u_1$ to the path $P_T(v_1,v_k)$ and $P_T(w_1,v_{i_2})$ is the path in $T$ from $w_1$
to the path $P_T(v_1,v_k)$. Then when $j\geq \max\{i_1,i_2\}$, $u_1$ and $w_1$ are
contained in $S_j$. We take the edge $v_jv_{j+1}$ in $P_T(v_1,v_k)$ with such index $j$.
Since $u_1,w_1\in S_j$, it follows that there exists an edge $u_hu_{h+1}$ of $P^A_{n_d}$
which is contained in $\partial (S_j)$ and there exists an edge $w_lw_{l+1}$ of $P^B_{n_d}$
which is contained in $\partial (S_j)$. Therefore, these edges $u_hu_{h+1}$ and $w_lw_{l+1}$
are cotree edges in $\overline T$.

On the other hand, for the tree-edge $v_jv_{j+1}\in T$, we may assume that $v_j\in V(H_a)$ and
$v_{j+1}\in V(H_{a+1})$. For otherwise ($v_j$ and $v_{j+1}$ belong to the same $H_a$) we can
take a greater $j$. Then these two vertices can be represented by $v_j=(x_1,x_2,\ldots x_{d-1},a)$
and $v_{j+1}=(x_1,x_2,\ldots x_{d-1},a+1)$. In the component $T[S_j]$ of $T-v_jv_{j+1}$, we have
two paths $P_T(u_h,v_j)$ and $P_T(w_l,v_j)$, while in the component $T[\overline S_j]$ of
$T-v_jv_{j+1}$, we have two paths $P_T(u_{h+1},v_{j+1})$ and $P_T(w_{l+1},v_{j+1})$. In this
way, we obtain two fundamental cycles
$$P_T(u_h,v_j)\cup \{v_jv_{j+1}\}\cup P_T(u_{h+1},v_{j+1})\cup \{u_hu_{h+1}\}$$
and
$$P_T(w_l,v_j)\cup \{v_jv_{j+1}\}\cup P_T(w_{l+1},v_{j+1})\cup \{w_lw_{l+1}\}$$
where $u_hu_{h+1}$ and $w_lw_{l+1}$ are cotree edges. The maximum stretch incurred by these 
two fundamental cycles is at least
$$\max\{2\sum_{i=1}^{d-1}(x_i-1)+1, 2\sum_{i=1}^{d-1}(n_i-x_i)+1\}\geq
2\sum_{i=1}^{d-1}\left\lceil\frac{n_i-1}{2}\right\rceil+1$$
$$=2\left(\left\lfloor\frac{n_1}{2}\right\rfloor+\left\lfloor\frac{n_2}{2}\right\rfloor
+\cdots+\left\lfloor\frac{n_{d-1}}{2}\right\rfloor\right)+1.$$
This yields the required lower bound.

{\bf Case 2:} \ $P^A_{n_d}$ is contained in $T$ but $P^B_{n_d}$ is not.

We use $P^A_{n_d}$ in place of $P_T(v_1,v_k)$. For an edge $u_ju_{j+1}\in T$ of $P^A_{n_d}$,
let $S_j$ and $\overline S_j$ be the vertex sets of the two components of $T-u_ju_{j+1}$,
where $u_j\in S_j$ and $u_{j+1}\in \overline S_j$. Then $\partial (S_j)$ is a fundamental
edge-cut with respect to $u_ju_{j+1}\in T$. We can choose such $j$ that $w_1=(n_1,n_2,\ldots,
n_{d-1},1)$ is contained is $S_j$. Hence there exists an edge $w_lw_{l+1}$ of $P^B_{n_d}$
which is contained in $\partial (S_j)$, and so it is a cotree edge in $\overline T$. On the
other hand, suppose that $u_j=(1,1,\ldots 1,a)$ and $u_{j+1}=(1,1,\ldots 1,a+1)$. Then we can
find a path $P_T(w_l,u_j)$ in the component $T[S_j]$ of $T-u_ju_{j+1}$ and a path $P_T(w_{l+1},
u_{j+1})$ in the component $T[\overline S_j]$. This results in a fundamental cycle
$$P_T(w_l,u_j)\cup \{u_ju_{j+1}\}\cup P_T(w_{l+1},u_{j+1})\cup \{w_lw_{l+1}\}.$$
Consequently, the stretch incurred by this fundamental cycle is at least
$$2\sum_{i=1}^{d-1}(n_i-1)+1>2(\left\lfloor\frac{n_1}{2}\right\rfloor+\left\lfloor\frac{n_2}{2}\right\rfloor
+\cdots+\left\lfloor\frac{n_{d-1}}{2}\right\rfloor)+1,$$
as required.

{\bf Case 3:} \ Both $P^A_{n_d}$ and $P^B_{n_d}$ are contained in $T$.

If all pairs of antipodal boundary paths are contained is $T$, then there would be cycles in $T$,
which contradicts that $T$ is a spanning tree. Otherwise we can take a pair antipodal boundary
paths $P^A_{n_i}$ and $P^B_{n_i}$ in some $x_i$-coordinate, not both of which are contained in
$T$. Now suppose that
\begin{eqnarray*}
P^A_{n_i} &:=& ((1,\ldots,1,1,1,\ldots,1),(1,\ldots,1,2,1\ldots,1),\ldots,(1,\ldots,1,n_i,1\ldots,1)),\\
P^B_{n_i} &:=& ((n_1,\ldots,n_{i-1},1,n_{i+1}\ldots,n_d),(n_1,\ldots,n_{i-1},2,n_{i+1}\ldots,n_d),\\
&& \ldots,(n_1,\ldots,n_{i-1},n_i,n_{i+1}\ldots,n_d)).
\end{eqnarray*}
By the same method of Case 1 and Case 2 for $P^A_{n_i}$ and $P^B_{n_i}$, we can show that
$$\sigma_T(G)\geq 2\sum_{1\leq j\leq d,j\neq i}\left\lfloor\frac{n_j}{2}\right\rfloor+1$$
$$\geq 2(\left\lfloor\frac{n_1}{2}\right\rfloor+\left\lfloor\frac{n_2}{2}\right\rfloor
+\cdots+\left\lfloor\frac{n_{d-1}}{2}\right\rfloor)+1,$$
since $n_i\leq n_d$. This completes the proof. $\Box$

\section{Concluding remarks}
\hspace*{0.3cm} The minimum stretch spanning tree problem and the minimum congestion
spanning tree problem are two dual problems in spanning tree optimization, in which the
fundamental cycles and the fundamental edge-cuts are considered respectively. They have
applications in information science and have close relations to labeling and embedding
for graphs (such as the bandwidth and cutwidth problem). It is meaningful to establish
connections between these two problems. What we have seen from the above are the exact
formulas of $\sigma_T(G)$ for two families of graphs, the Hamming graphs and
higher-dimensional grids. The corresponding results for the tree-congestion $c_T(G)$
have not been seen yet.

The study of these optimization problems gives rise to two graph-theoretic invariants,
the tree-stretch $\sigma_T(G)$ and tree-congestion $c_T(G)$. From the perspective of
graph theory, several aspects are worthwhile to explored. For example, exact
representations for more graph families, relations with other parameters, extremal
graph characterizations, duality, symmetry, decomposability, etc., are expected.


\begin{thebibliography}{abc}

\bibitem{Bondy08} J.A. Bondy and U.S.R. Murty, Graph Theory.
Springer-Verlag, Berlin, 2008.

\bibitem{Bondy89}  J.A. Bondy, Trigraphs. {\it Discrete Math.},
75(1-3)(1989) 69-79.

\bibitem{Bodlae11} H.L. Bodlaender, K. Kozawa, T. Matsushima and
Y. Otachi, Spanning tree congestion of $k$-outerplanar graphs. {\it
Discrete Math.}, 311 (2011) 1040-1045.

\bibitem{Bodlae12} H.L. Bodlaender, F.V. Fomin, P.A. Golovach,
Y.Otachi and  E.J. van Leeuwen, Parameterized complexity of the
spanning tree congestion problem. {\it Algorithmica}, 64 (2012)
85-111.

\bibitem{Brand04} A. Brandst\"adt, F.F. Dragan, H.-O. Le and
V.B. Le, Tree spanners on chordal graphs: complxity and algorithms.
{\it Theor. Comput. Sci.}, 310 (2004) 329-354.

\bibitem{Brand07}  A. Brandst\"adt, F.F. Dragan, H.-O. Le,
V.B. Le and  R. Uehara, Tree spanners for bipartite graphs and probe
interval graphs. {\it Algorithmica}, 47 (2007) 27-51.

\bibitem{Cai95}  L. Cai and D.G. Corneil, Tree spanners. {\it SIAM
J. Discrete Math.}, 8 (1995) 359-387.

\bibitem{Caste09}  A. Castej\'on and  M.I. Ostrovskii, Minimum congestion
spanning trees of grids and discrete toruses. {\it Discuss. Math.
Graph Theory}, 29(2009) 511-519.

\bibitem{Diaz02} J. Diaz, J. Petit and  M. Serna, A survey of
graph layout problems. {\it ACM Computing Surveys}, 34 (2002)
313-356.

\bibitem{Fekete01} S.P. Fekete and J. Kremer, Tree spanners in
planar graphs. {\it Discrete Appl. Math.}, 108 (2001) 85-103.

\bibitem{Fomin11} F.V. Fomin, P.A. Golovach and E.J. van Leeuwen,
Spanners of bounded degree graphs. {\it Inform. Process. Lett.}, 111
(2011) 142-144.

\bibitem{Galb03} G. Galbiati, On finding cycle bases and fundamental
cycle bases with a shortest maximal cycles. {\it Inform. Process.
Lett.}, 88 (2003) 155-159.

\bibitem{Hruska08} S.W. Hruska, On tree congestion of graphs.
{\it Discrete Math.},  308 (2008) 1801-1809.

\bibitem{Le99} H.-O. Le and V.B. Le, Optimal tree 3-spanners in
directed path graphs. {\it Networks}, 34 (1999) 81-87.

\bibitem{Lieb08}  C. Liebchen and G. W\"unsch, The zoo of tree
spanner problems. {\it Discrete Appl. Math.}, 156 (2008) 569-587.

\bibitem{Lin17} L. Lin and Y. Lin, The minimum stretch spanning
tree problem for typical graphs. To appear in Acta Math. Appl. Sinica
(arXiv:1712.03497v1 [math.CO] 10 Dec 2017).

\bibitem{Lin18} L. Lin and Y. Lin, The minimum stretch spanning
tree problem for graph products. Summitted 2018.

\bibitem{Madan96} M.S. Madanlal, G. Venkatesan and  C.P. Rangan,
Tree 3-apanners on interval, permutation and regular bipartite
graphs. {\it Inform. Process. Lett.}, 59 (1996) 97-102.

\bibitem{Okamoto11} Y. Okamoto, Y. Otachi, R. Uehara, and T. Uno,
Hardness results and an exact exponential algorithm for the spanning
tree congestion problem, {\it J. of Graph Algorithms and
Applications}, 15(6)(2011) 727-751.

\bibitem{Ostrov04}  M.I. Ostrovskii, Minimal congestion trees.
{\it Discrete Math.},  285 (2004) 219-326.

\bibitem{Ostrov10} M.I. Ostrovskii, Minimum congestion spanning trees
in planar graphs. {\it Discrete Math.}, 310 (2010) 1204-1209.

\bibitem{Peleg89} D. Peleg and J.D. Ullman, An optimal synchroniser
for the hypercube. {\it SIAM J. Comput.}, 18(4) (1989) 740-747.

\bibitem{Venka97}  G. Venkatesan, U. Rotics, M.S. Madanlal, J.A. Makowsy,
and  C.P. Rangan, Restrictions of minimum spanner problems. {\it
Information and Computation}, 136(1997) 143-164.



\end{thebibliography}
\end{document}